\documentclass[11pt]{article}
\usepackage{pifont}

\usepackage{amsfonts}
\usepackage{latexsym}
\usepackage{amsmath}
\usepackage{amssymb}
\usepackage{color}

\newcommand{\BR}{\mathbb{R}}

\newcommand{\BN}{\mathbb{N}}
\newcommand{\CN}{\mathcal{N}}

\newcommand{\RN}{\mathbb{R}^n}

\newcommand{\E}{\hat{\mathbb{E}}}
\newcommand{\CH}{\mathcal{H}}
\newcommand{\CLL}{C_{l,lip}(\mathbb{R}^n)}

\newcommand{\FI}{\varphi}
\newcommand{\GX}{\underline{\sigma}^2}
\newcommand{\GS}{\bar{\sigma}^2}

\newtheorem{theorem}{Theorem}[section]

\newtheorem{proposition}{Proposition}[section]

\newtheorem{corol}{Corollary}[section]

\newtheorem{definition}{Definition}[section]
\newtheorem{remark}{Remark}[section]

\title{Multiple G-It\^{o} integral in the G-expectation space}

\date{}

\author{Panyu Wu\thanks{E-mail:panyuaza@yahoo.com.cn}\\Department of Mathematics,\\Shandong University, Jinan, China}

\begin{document}

\maketitle

\begin{abstract}
In this paper, motivated by mathematic finance we introduce the multiple G-It\^{o} integral in the G-expectation space, then investigate how to calculate. We get the the relationship between Hermite polynomials and multiple G-It\^{o} integrals which is a natural extension of the classical result obtained by It\^{o} in 1951.
\par  $\textit{Keywords:}$ Sublinear expectation, G-Brownian motion, G-It\^{o} integral, Hermite polynomials
\end{abstract}



\section{\textbf{Introduction}}
A multiple stochastic integral with respect to the classical Brownian motion was constructed by Wiener in Ref. \cite{wie} as a polynomial chaos in independent
Gaussian random variables. A more general construction was due to It\^{o} in Ref. \cite{ito}. Actually, the theory and applications of It\^{o} multiple stochastic
integrals are fairly rich, for example, Engel \cite{en} for the history and framework of multiple integration, Cheridito et al. \cite{che} for applications
in finance and Soner et al. \cite{soner} for applications in stochastic target problems. However, the classical Brownian motion was constructed in a
linear expectation space, such linearity assumption is not feasible in many areas of applications because many uncertain phenomena can not be well modelled
using additive probabilities or linear expectations. More specifically, motivated by the risk measures and stochastic volatility problems in finance,
Peng in Ref. \cite{p2006} introduced the sublinear expectation space and initiated the G-normal distribution under a sublinear expectation space.
He also introduced the notions of G-Brownian motion as the counterpart of classical Brownian motion in the linear case and G-It\^{o} integral with respect to G-Brownian motion. He introduced a class of sublinear expectation space called G-expectation space as well and proved there exist G-Brownian motion in G-expectation space.
Now more and more people are interested in G-expectation space and the applications of G-It\^{o} integral will be more and more widely. A natural question is the following: how to define and calculate the multiple G-It\^{o} integral. The purpose of this paper is to solve this problem. We not only introduce the multiple G-It\^{o} integral of symmetric function in $L^2([0,T]^n)$ but also obtain the relationship between Hermite polynomials and multiple G-It\^{o} integrals. All of them are natural and fairly neat extensions of the classical It\^{o}'s results, but the proof here is different from the original proof of the classical multiple It\^{o} integrals.
\par
The remainder of this paper is organized as follows. In section 2, we recall some notions and results in G-expectation space which will be useful in this paper.
In section 3, we introduce the multiple G-It\^{o} integral. In section 4, we state and prove the main result of this paper which is the relationship between
Hermite polynomials and multiple G-It\^{o} integrals.
\section{\textbf{Preliminaries}}
We recall some notions and results in G-expectation space. Some more details can be found in Refs. [4-8].
\begin{definition}
A random variable $X$ on a sublinear expectation space $(\Omega,\CH,\E)$ is called G-normal distributed, denoted by $X\sim \CN(0,[\underline{\sigma}^2,\GS])$, if $$aX+b\bar{X}\sim\sqrt{a^2+b^2}X,\ \forall a,b\geq0,$$ where $\bar{X}$ is an independent copy of $X$, $\GS=\E[X^2]$ and $\underline{\sigma}^2=-\E[-X^2].$ Here the letter G denotes the function $G(\alpha):={1\over 2}\E[\alpha X^2]={1\over2}(\GS\alpha^+-\GX\alpha^-):\BR\rightarrow\BR.$
\end{definition}
\begin{definition}
Let $G(\cdot):\BR\rightarrow\BR,G(\alpha)={1\over2}(\GS\alpha^+-\GX\alpha^-)$, where $0\leq\underline{\sigma}\leq\bar{\sigma}<\infty.$ A stochastic process $(B_t)_{t\geq0}$ in a sublinear expectation space $(\Omega,\CH,\E)$
is called a G-Brownian motion if the following properties are satisfied:\\
(i) $B_0(\omega)=0;$\\
(ii) For each $t,s\geq0$, the increment $B_{t+s}-B_t$ is $\CN(0,[s\underline{\sigma}^2,s\GS])$-distributed
and is independent to $(B_{t_1},B_{t_2},\cdots,B_{t_n})$, for each $n\in\BN$ and $0\leq t_1\leq t_2\leq\cdots\leq t_n\leq t.$
\end{definition}
In the rest of this paper, we denote by $\Omega=C_0(\BR^+)$ the space of all $\BR$-valued continuous paths $(\omega_t)_{t\in\BR^+}$ with $\omega_0=0$, equipped with the distance
$$\rho(\omega^1,\omega^2):=\sum_{i=1}^{\infty}2^{-i}[(\max_{t\in[0,i]}|\omega_t^1-\omega_t^2|)\wedge1].$$
For each fixed $T\in[0,\infty),$ we set $\Omega_T:=\{\omega_{\cdot\wedge T}:\omega\in\Omega\},$
$$L_{ip}(\Omega_T):=\{\FI(B_{t_1\wedge T},\cdots,B_{t_n\wedge T}):n\in \BN,t_1,\cdots,t_n\in[0,\infty),\FI\in\CLL\},$$
$L_{ip}(\Omega):=\cup_{n=1}^{\infty}L_{ip}(\Omega_n),$
where $B_t$ denote the canonical process, that is, $B_t(\omega)=\omega_t.$
\par
For any given monotonic and sublinear function $G(\cdot):\BR\rightarrow\BR$, consider the G-expectation $\E[\cdot]:L_{ip}(\Omega)\rightarrow\BR$ defined by Peng in Ref. \cite{p2006}. He proved that the corresponding canonical process $(B_t)_{t\geq0}$ on the sublinear expectation space $(\Omega,L_{ip}(\Omega),\E)$ called G-expectation space is a G-Brownian motion. In the sequel, G-Brownian motion means the canonical process $(B_t)_{t\geq0}$ under the G-expectation $\E$.
\par
We denote the completion of $L_{ip}(\Omega)$ under the norm $\| X\| _p:=(\E[|X|^p])^{{1\over p}}$ by $L_G^p(\Omega),p\geq 1.$ Let $M_G^{p,0}(0,T)$ be the collection of processes in the following form:
$$\eta_t(\omega)=\sum_{k=0}^{N-1}\xi_k(\omega)I_{[t_k,t_{k+1})}(t),$$
where ${0=t_0<t_1<\cdots< t_N= T}$ is any given partition  of $[0,T],\xi_k\in L_G^p(\Omega_{t_k}),k=0,\cdots,N-1$.
For each $\eta\in M_G^{p,0}(0,T)$, let $\|\eta\|_{M_G^p}=[\E(\int_0^T|\eta_s|^pds)]^{{1\over p}}$ and $M_G^p(0,T)$ denote the completion of $M_G^{p,0}(0,T)$ under
norm $\|\cdot\|_{M_G^p}$.
\par
Let $(B_t)_{t\geq0}$ be a G-Brownian motion with $G(\alpha)={1\over2}(\GS\alpha^+-\GX\alpha^-)$, where $0\leq\underline{\sigma}\leq\bar{\sigma}<\infty.$
\begin{definition}
For each $\eta\in M_G^{2,0}(0,T)$ of the form $\eta_t(\omega)=\sum_{k=0}^{N-1}\xi_k(\omega)I_{[t_k,t_{k+1})}(t)$, we define
$$I(\eta)=\int_0^T\eta_tdB_t:=\sum_{j=0}^{N-1}\xi_j(B_{t_{j+1}}-B_{t_j}).$$
\end{definition}
\begin{proposition}
The mapping $I(\cdot):M_G^{2,0}(0,T)\rightarrow L_G^2(\Omega_T)$ is a continuous linear mapping under norm $\|\cdot\|_{M_G^2}$ and $\|\cdot\|_2$, thus $I(\cdot)$ can be continuously extended to $M_G^2(0,T)$. For any $\eta\in M_G^2(0,T)$, we denote $\int_0^T\eta_tdB_t:=I(\eta).$ And we have
\begin{equation}\label{eq0}
\E[(\int_0^T\eta_tdB_t)^2]\leq\GS\E[\int_0^T\eta_t^2dt].
\end{equation}
\end{proposition}
\begin{definition}
The quadratic variation process of G-Brownian motion $(B_t)_{t\geq0}$ is defined by
$$\langle B\rangle_t:=\lim_{\mu(\pi_t^N)\rightarrow 0}\sum_{j=0}^{N-1}(B_{t_{j+1}^N}-B_{t_j^N})^2,$$
where $\mu(\pi_t^N):=\max\{|t_{i+1}^N-t_i^N|:0=t_0<t_1<\cdots< t_N= t\}.$
\end{definition}

\begin{definition}
For each $\eta\in M_G^{1,0}(0,T)$ of the form $\eta_t(\omega)=\sum_{k=0}^{N-1}\xi_k(\omega)I_{[t_k,t_{k+1})}(t)$, we define
$$Q(\eta)=\int_0^T\eta_td\langle B\rangle_t:=\sum_{j=0}^{N-1}\xi_j(\langle B\rangle_{t_{j+1}}-\langle B\rangle_{t_j}):M_G^{1,0}(0,T)\rightarrow L_G^1(\Omega_T).$$
\end{definition}
\begin{proposition}
The mapping $Q(\cdot):M_G^{1,0}(0,T)\rightarrow L_G^1(\Omega_T)$ is a continuous linear mapping under norm $\|\cdot\|_{M_G^1}$ and $\|\cdot\|_1$, thus $Q(\cdot)$ can be continuously extended to $M_G^1(0,T)$. For any $\eta\in M_G^2(0,T)$, we denote $\int_0^T\eta_td\langle B\rangle_t:=Q(\eta).$
\end{proposition}

\begin{proposition}\textbf{G-It\^{o}'s formula:} Let $\Phi\in C^2(\RN)$ with $\partial^2_{x_ix_j}\Phi$ satisfying polynomial growth condition for $i,j=1,\cdots,n$, and $X_t=(X_t^1,\cdots,X_t^n)$ satisfying
$$X_t^i=X_0^i+\int_0^t\alpha_s^ids+\int_0^t\eta_s^i d\langle B\rangle_s +\int_0^t \beta_s^i dB_s,\ \ i=1,\cdots,n,$$
where $\alpha^i,\eta^i,\beta^i$ be bounded processes in $M_G^2(0,T).$ Then for each $t\geq 0$ we have, in $L_G^2(\Omega_t)$
\begin{eqnarray*}
\Phi(X_t)-\Phi(X_s)&=&\sum_{i=1}^n[\int_s^t\partial_{x_i}\Phi(X_u)\alpha_u^idu+\int_s^t\partial_{x_i}\Phi(X_u)\beta_u^i dB_u]\\
& &+\int_s^t[\sum_{i=1}^n\partial_{x_i}\Phi(X_u)\eta_u^i+{1\over2}\sum_{i,j=1}^n\partial^2_{x_ix_j}\Phi(X_u)\beta_u^i\beta_u^j]d\langle B\rangle_u.
\end{eqnarray*}
\end{proposition}
\section{\bf{Multiple G-It\^{o} integrals}}
In order to introduce the definition of multiple G-It\^{o} integral, we introduce the following usual spaces of function:\\
$L^2([0,T]^n):=\{g|g:[0,T]^n\rightarrow\BR,\|g\|^2_{L^2([0,T]^n)}<\infty\};$\\
$\hat{L}^2([0,T]^n):=\{g|g$ is a symmetric function in $L^2([0,T]^n) \},$\\
where $\|g\|^2_{L^2([0,T]^n)}=\int_{[0,T]^n}g^2(x_1,\cdots,x_n)dx_1\cdots dx_n.$
\par
For any $f$ on $S_n:=\{(x_1,\cdots,x_n)\in[0,T]^n:0\leq x_1\leq x_2\leq \cdots\leq x_n\leq T\}\ (n\geq 1)$, we define
$$\|f\|^2_{L^2(S_n)}:=\int_{S_n}f^2(x_1,\cdots,x_n)dx_1\cdots x_n.$$
\par
For $\|f\|^2_{L^2(S_n)}<\infty$ we can form the ($n-$fold) iterated G-It\^{o} integral
$$J_n^T(f):=\int_0^T\int_0^{t_n}\cdots\int_0^{t_3}\int_0^{t_2}f(t_1,t_2,\cdots,t_n)dB_{t_1}dB_{t_2}\cdots dB_{t_n}.$$
It is easily to show that at each G-It\^{o} integration with respect to $dB(t_i)$ is included in $M_G^2(0,t_{i+1})$ by equality (\ref{eq0}). Moreover, by equality (\ref{eq0}) we have
\begin{eqnarray*}
\E[(J_n^T(f))^2]&=&\E\left[\left(\int_0^T\int_0^{t_n}\cdots\int_0^{t_2}f(t_1,t_2,\cdots,t_n)dB_{t_1}dB_{t_2}\cdots dB_{t_n}\right)^2\right]\\
&\leq&\GS\int_0^T\E\left[\left(\int_0^{t_n}\cdots\int_0^{t_2}f(t_1,t_2,\cdots,t_n)dB_{t_1}dB_{t_2}\cdots dB_{t_{n-1}}\right)^2\right]dt_n\\
&\leq&\bar{\sigma}^{2n}\int_0^T\int_0^{t_n}\cdots\int_0^{t_2}f^2(t_1,t_2,\cdots,t_n)dt_1\cdots dt_n\\
&=&\bar{\sigma}^{2n}\|f\|^2_{L^2(S_n)}<\infty.
\end{eqnarray*}
\par
For any constant $c$, we define $J_0(c)=c$. Notice that for any $g\in \hat{L}^2([0,T]^n)$, we have
$$\|g\|^2_{L^2(S_n)}={1\over{n!}}\|g\|^2_{L^2([0,T]^n)}.$$
\par
Thus we give the following definition of multiple G-It\^{o} integral.
\begin{definition}\label{de1}
For any $g\in \hat{L}^2([0,T]^n)$, define
$$I_n^T(g):=\int_{[0,T]^n}g(t_1,\cdots,t_n)dB_{t_1}dB_{t_2}\cdots dB_{t_n}:=n!J_n^T(g)$$
\end{definition}
Notice that for all $g\in \hat{L}^2([0,T]^n)$, we have $I_n^T(g)\in L^2_G(\Omega_T)$ because of
$$\E[(I^T_n(g))^2]=\E[(n!)^2(J^T_n(g))^2]\leq\bar{\sigma}^{2n}(n!)^2\|g\|^2_{L^2(S_n)}=\bar{\sigma}^{2n}n!\|g\|^2_{L^2([0,T]^n)}.$$

\section{\bf{Main Result}}
We start by introducing the Hermite polynomials $h_n(x)$ which are defined by
\begin{equation}\label{eq1}
h_n(x)=(-1)^ne^{{1\over2}x^2}{d^n\over{dx^n}}(e^{-{1\over2}x^2}),\ \ n=0,1,2,\cdots.
\end{equation}
Obviously the first three Hermite polynomials are:
$$h_0(x)=1,\ h_1(x)=x,\ h_2(x)=x^2-1.$$
\par
We claim the main result as the following theorem:
\begin{theorem}\label{th1}
For any $f\in L^2([0,T])$, let $g_n(t_1,t_2,\cdots,t_n)=f(t_1)f(t_2)\cdots f(t_n)$, then $g_n\in \hat{L}^2([0,T]^n)$, and in $L^2_G(\Omega_T)$
\begin{equation}\label{eq2}
I_n^T(g_n)=\|f\|_T^nh_n({\theta_T\over{\|f\|_T}}),
\end{equation}
where $\|f\|_T=[\int_0^Tf^2(s)d\langle B\rangle_s]^{1\over2}$ be a nonnegative random variable and $\theta_T=\int_0^Tf(t)dB_t.$
\end{theorem}
\hspace{-0.7cm}{\bf{\it Proof}}\hspace{0.4cm}
It is easy to check that $g\in \hat{L}^2([0,T]^n)$. We now prove the theorem in two steps.
\par
Step 1: Equality (\ref{eq2}) holds if and only if the following equality (\ref{eq3}) is true:
\begin{equation}\label{eq3}
I_n^T(g_n)=\theta_T I^T_{n-1}(g_{n-1})-(n-1)\|f\|_T^2I^T_{n-2}(g_{n-2}),\ \ n\geq2.
\end{equation}
\par
On the one hand, if equality (\ref{eq2}) holds, using the Hermite polynomials's recurrence relation:
\begin{equation}\label{eq4}
h_n(y)=y h_{n-1}(y)-(n-1)h_{n-2}(y),\ \ n\geq2.
\end{equation}
For $n\geq2$ we have:
\begin{eqnarray*}
I_n^T(g_n)&=&\|f\|_T^nh_n({\theta_T\over{\|f\|_T}})\\
&=&\theta_T \|f\|_T^{n-1}h_{n-1}({\theta_T\over{\|f\|_T}})-(n-1)\|f\|_T^nh_{n-2}({\theta_T\over{\|f\|_T}})\\
&=&\theta_T I^T_{n-1}(g_{n-1})-(n-1)\|f\|_T^2I^T_{n-2}(g_{n-2}).
\end{eqnarray*}
We obtain that the equality (\ref{eq3}) holds for any $n\geq2$.
\par
On the other hand, if equality (\ref{eq3}) holds, obviously we have
\begin{equation}\label{eq5}
I_0^T(g_0)=1=\|f\|^0h_0({\theta_T\over{\|f\|_T}}),
\end{equation}
\begin{equation}\label{eq6}
I_1^T(g_1)=\theta_T=\|f\|_Th_1({\theta_T\over{\|f\|_T}}).
\end{equation}
When $n=2$, applying G-It\^{o}'s formula to $\theta_t^2$, we get $$d\theta_t^2=2(\int_0^tf(s)dB_s)f(t)dB_t+f^2(t)d\langle B\rangle_t,$$
that is $\theta_T^2=2\int_0^T\int_0^{t}f(s)f(t)dB_sdB_t+\int_0^Tf^2(t)d\langle B\rangle_t.$
Hence,
\begin{eqnarray*}
\|f\|_T^2h_2({\theta_T\over{\|f\|_T}})&=&\|f\|_T^2[({\theta_T\over{\|f\|_T}})^2-1]\\
&=&\theta_T^2-\int_0^Tf^2(t)d\langle B\rangle_t\\
&=&2\int_0^T\int_0^{t}f(s)f(t)dB_sdB_t.
\end{eqnarray*}
Therefore,
\begin{equation}\label{eq7}
\|f\|_T^2h_2({\theta_T\over{\|f\|_T}})=I^T_2(g_2).
\end{equation}
From equality (\ref{eq5})-(\ref{eq7}) it follows that equality (\ref{eq2}) holds for $n=0,1,2$. For $n>2$ equality (\ref{eq2}) can easily been proved by mathematical induction using equality (\ref{eq3}) and (\ref{eq4}), we omit it.
\par
Step 2: We shall show that equality (\ref{eq3}) holds under the assumption of the theorem. We deduce from equations (\ref{eq5})-(\ref{eq7}) that equation (\ref{eq3}) holds true in case $n=2$. We make use of the mathematical induction with regard $n$. Now suppose equality (\ref{eq3}) holds when $n\leq m-1$, we have to prove equality (\ref{eq3}) being true when $n=m$.
\par
Let
$$X_t=\int_0^t\int_0^{t_{m-1}}\cdots\int_0^{t_2}f(t_1)\cdots f(t_{m-1})dB_{t_1}\cdots dB_{t_{m-1}}.$$
By G-It\^{o}'s formula, we get
\begin{eqnarray*}
d\theta_t X_t&=&\left(\int_0^t\int_0^{t_{m-1}}\cdots\int_0^{t_2}f(t_1)\cdots f(t_{m-1})dB_{t_1}\cdots dB_{t_{m-1}}\right)f(t)dB_t\\
& & +\theta_t\left(\int_0^t\int_0^{t_{m-2}}\cdots\int_0^{t_2}f(t_1)\cdots f(t_{m-2})dB_{t_1}\cdots dB_{t_{m-2}}\right)f(t)dB_t\\
& & +\left(\int_0^t\int_0^{t_{m-2}}\cdots\int_0^{t_2}f(t_1)\cdots f(t_{m-2})dB_{t_1}\cdots dB_{t_{m-2}}\right)f^2(t)d\langle B\rangle_t.
\end{eqnarray*}
Thus,
\begin{eqnarray*}
&&\theta_T I^T_{m-1}(g_{m-1})\\
&&\hspace*{-1em}=(m-1)!\theta_T X_T\\
&&\hspace*{-1em}=(m-1)!\int_0^T \int_0^t\int_0^{t_{m-1}}\cdots\int_0^{t_2}f(t_1)\cdots f(t_{m-1})f(t)dB_{t_1}\cdots dB_{t_{m-1}}dB_t\\
&&+(m-1)!\int_0^T \left(\int_0^t\int_0^{t_{m-2}}\cdots\int_0^{t_2}f(t_1)\cdots f(t_{m-2})dB_{t_1}\cdots dB_{t_{m-2}}\right)f^2(t)d\langle B\rangle_t\\
&&+(m-1)!\int_0^T\theta_t\left(\int_0^t\int_0^{t_{m-2}}\cdots\int_0^{t_2}f(t_1)\cdots f(t_{m-2})dB_{t_1}\cdots dB_{t_{m-2}}\right)f(t)dB_t\\
&&\hspace*{-1em}=I_m^T(g_m)+(m-1)!\int_0^T \left(\int_0^t\int_0^{t_{m-2}}\cdots\int_0^{t_2}f(t_1)\cdots f(t_{m-2})dB_{t_1}\cdots dB_{t_{m-2}}\right)f^2(t)d\langle B\rangle_t\\
&&+(m-1)!\int_0^T\theta_t\left(\int_0^t\int_0^{t_{m-2}}\cdots\int_0^{t_2}f(t_1)\cdots f(t_{m-2})dB_{t_1}\cdots dB_{t_{m-2}}\right)f(t)dB_t\\
&&-(m-1)(m-1)!\int_0^T \int_0^t\int_0^{t_{m-1}}\cdots\int_0^{t_2}f(t_1)\cdots f(t_{m-1})f(t)dB_{t_1}\cdots dB_{t_{m-1}}dB_t\\
&&\hspace*{-1em}=I_m^T(g_m)+(m-1)!\int_0^T \left(\int_0^t\int_0^{t_{m-2}}\cdots\int_0^{t_2}f(t_1)\cdots f(t_{m-2})dB_{t_1}\cdots dB_{t_{m-2}}\right)f^2(t)d\langle  B\rangle_t\\
&&+(m-1)\psi_m,
\end{eqnarray*}
where
\begin{eqnarray*}
&&\psi_m=\int_0^T\left[(m-2)!\theta_t\left(\int_0^t\int_0^{t_{m-2}}\cdots\int_0^{t_2}f(t_1)\cdots f(t_{m-2})dB_{t_1}\cdots dB_{t_{m-2}}\right)\right.\\
&&\hspace*{3em}\left.-(m-1)!\int_0^t\int_0^{t_{m-1}}\cdots\int_0^{t_2}f(t_1)\cdots f(t_{m-1})dB_{t_1}\cdots dB_{t_{m-1}}\right]f(t)dB_t\\
&&\hspace*{1.8em}=\int_0^T[\theta_tI^t_{m-2}(g_{m-2})-I^t_{m-1}(g_{m-1})]f(t)dB_t.
\end{eqnarray*}
From equality (\ref{eq3}), we have
\begin{eqnarray*}
&&\psi_m=\int_0^T[(m-2)\|f\|^2_tI^t_{m-3}(g_{m-3}))]f(t)dB_t\\
&&\hspace*{-1em}=(m-2)!\int_0^T\left(\int_0^tf^2(s)d\langle B\rangle_s\right)\left(\int_0^t\int_0^{t_{m-3}}\cdots\int_0^{t_2}f(t_1)\cdots f(t_{m-3})dB_{t_1}\cdots dB_{t_{m-3}}\right)f(t)dB_t.
\end{eqnarray*}
\par
Let $Y_t=\int_0^t\int_0^{t_{m-2}}\cdots\int_0^{t_2}f(t_1)\cdots f(t_{m-2})dB_{t_1}\cdots dB_{t_{m-2}}$, using G-It\^{o}'s formula to $\|f\|^2_tY_t$, we get:
\begin{eqnarray*}
&&(m-1)\|f\|_T^2I^T_{m-2}(g_{m-2})=(m-1)!\|f\|^2_TY_T\\
&&\hspace*{-1em}=(m-1)!\int_0^T \left(\int_0^t\int_0^{t_{m-2}}\cdots\int_0^{t_2}f(t_1)\cdots f(t_{m-2})dB_{t_1}\cdots dB_{t_{m-2}}\right)f^2(t)d\langle  B\rangle_t\\
&&\hspace*{-1em}+(m-1)!\int_0^T\left(\int_0^tf^2(s)d\langle B\rangle_s\right)\left(\int_0^t\int_0^{t_{m-3}}\cdots\int_0^{t_2}f(t_1)\cdots f(t_{m-3})dB_{t_1}\cdots dB_{t_{m-3}}\right)f(t)dB_t\\
&&\hspace*{-1em}=(m-1)!\int_0^T \left(\int_0^t\int_0^{t_{m-2}}\cdots\int_0^{t_2}f(t_1)\cdots f(t_{m-2})dB_{t_1}\cdots dB_{t_{m-2}}\right)f^2(t)d\langle  B\rangle_t+(m-1)\psi_m.
\end{eqnarray*}
Therefore, $\theta_T I^T_{m-1}(g_{m-1})=I_m^T(g_m)+(m-1)\|f\|_T^2I^T_{m-2}(g_{m-2})$, in other words, the equality (\ref{eq3}) has been established for $n=m$. By
mathematical induction, equality (\ref{eq3}) holds for any integer $n\geq2$. The proof of theorem \ref{th1} is complete.  $\square$
\par
\begin{remark} G-Brownian motion degenerate to the classical Brownian motion when $\GS=\GX=1.$ In that case, equality (\ref{eq2}) becomes the relation between the classical  multiple It\^{o} integrals and Hermite polynomials.
\end{remark}
\par
The next corollary gives the general formula of $\int_0^T\int_0^{t_n}\cdots\int_0^{t_2}dB_{t_1}\cdots dB_{t_n}$.
\begin{corol}\label{co}
$$\int_0^T\int_0^{t_n}\cdots\int_0^{t_2}dB_{t_1}\cdots dB_{t_n}=\sum_{m=0}^{\lfloor{n\over 2}\rfloor}{{(-1)^m}\over{2^m m!(n-2m)!}}\langle B\rangle_T^m B_T^{n-2m},$$
where $\lfloor x\rfloor$ is the largest integer not greater than $x$.
\end{corol}
\hspace{-0.7cm}{\bf{\it Proof}}\hspace{0.4cm}
From theorem \ref{th1} it follows that
$$\int_0^T\int_0^{t_n}\cdots\int_0^{t_2}dB_{t_1}\cdots dB_{t_n}={1\over{n!}}\langle B\rangle_T^{n\over2} h_n({B_T\over{\langle B\rangle_T^{1/2}}}).$$
It is easily to get the corollary since the Hermite polynomials can be written explicitly as
$$h_n(x)=n!\sum_{m=0}^{\lfloor{n\over 2}\rfloor}{{(-1)^m}\over{2^m m!(n-2m)!}}x^{n-2m}.$$
The proof of  corollary \ref{co} is complete.  $\square$

\end{document}